\documentclass{amsart} 
\usepackage{amscd}
\title{Homotopy invariance of AF-embeddability}
\author{Narutaka Ozawa} 
\address{Department of Mathematical Science,
University of Tokyo, 153-8914, Japan}
\email{narutaka@ms.u-tokyo.ac.jp} 
\date{January 20, 2002}
\newtheorem{thm}{Theorem} 
\newtheorem{prop}[thm]{Proposition}
\newtheorem{lem}[thm]{Lemma} 
\theoremstyle{definition} 
\newtheorem{defn}[thm]{Definition} 
\newcommand{\C}{{\mathbb C}} 
\newcommand{\N}{{\mathbb N}} 
\newcommand{\M}{{\mathbb M}} 
\newcommand{\B}{{\mathbb B}} 
\newcommand{\e}{\varepsilon} 
\newcommand{\p}{\varphi} 
\newcommand{\hh}{{\mathcal H}} 
\newcommand{\id}{\mathrm{id}} 
\newcommand{\diag}{\mathop{\mathrm{diag}}} 
\newcommand{\simm}{\mathop{\sim}}

\begin{document}
\begin{abstract}
We prove that AF-embeddability is a homotopy invariant 
in the class of separable exact $C^*$-algebras. 
This work was inspired by Spielberg's work on homotopy invariance of 
AF-embeddability and Dadarlat's serial works on AF-embeddability 
of residually finite dimensional $C^*$-algebras. 
\end{abstract}
\maketitle
\section{Introduction}
Kirchberg \cite{kcar} showed a remarkable theorem 
that a separable exact $C^*$-algebra is nuclearly embeddable and 
is $*$-isomorphic to a quotient of a subalgebra of the CAR algebra. 
(See \cite{dsbq} for a different proof.) 
However, the class of $C^*$-algebras 
which are $*$-isomorphic to subalgebras of AF-algebras 
(i.e., AF-embeddable $C^*$-algebras) is more restrictive. 
An AF-embeddable $C^*$-algebra has to be exact and quasidiagonal. 
Thus a major open problem on AF-embeddability is that 
whether the converse is true. 
This problem has the affirmative answer \cite{lin}, \cite{dene} 
for a separable residually finite dimensional (RFD) $C^*$-algebra 
which satisfies the universal coefficient theorem (UCT). 
We refer the reader to Wassermann's lecture notes \cite{wassermann} 
and Brown's survey \cite{brown} respectively for the information 
of exactness and of quasidiagonality. 
Spielberg \cite{spielberg} studied AF-embeddability of extensions 
and related AF-embeddability to homotopy types of $C^*$-algebras. 
Voiculescu \cite{voiculescu} then proved that quasidiagonality 
is a homotopy invariant. 
This, in turn, is used to show the homotopy invariance 
of AF-embeddability in the class of separable exact $C^*$-algebras. 
Following Voiculescu \cite{voiculescu}, we say a $C^*$-algebra $A$ 
is homotopically dominated by a $C^*$-algebra $B$ 
if there are $*$-homomorphisms $\p\colon A\to B$ and 
$\psi\colon B\to A$ such that $\psi\p$ is homotopic to $\id_A$. 
\begin{thm}
If $A$ is a separable exact $C^*$-algebra which is 
homotopically dominated by an AF-embeddable $C^*$-algebra $B$, 
then $A$ is AF-embeddable. 
\end{thm}

This theorem immediately follows from the next two propositions. 
The first is inspired by serial works \cite{drfd}, \cite{djfa}, 
\cite{dsbq}, \cite{dene} of Dadarlat on AF-embeddability of 
RFD $C^*$-algebras and the second is essentially same as Theorem 3.9 
of Spielberg \cite{spielberg} 
(namely, he proved it under an assumption on the UCT). 
We are indebted to Takeshi Katsura for clarifying Spielberg's 
argument to us. 

\begin{prop}\label{prop1}
The cone $CA$ over a separable exact $C^*$-algebra $A$ is AF-embeddable. 
\end{prop}

\begin{prop}\label{prop2}
Let $A$ be a separable $C^*$-algebra whose cone $CA$ is AF-embeddable. 
Suppose that $A$ is homotopically dominated by an AF-embeddable 
$C^*$-algebra $B$. 
Then $A$ is AF-embeddable. 
\end{prop}
\section{Proof of Proposition \ref{prop1}}
\begin{defn}
Let $F\subset C$ be a finite subset of a unital $C^*$-algebra $C$ 
and let $\e>0$. 
For a unital completely positive (ucp) map 
$\p$ from $C$ into a full matrix algebra $\M_p$, 
we give the following definitions. 

We say $\p$ is $(F,\e)$-multiplicative if 
$\|\p(f^*f)-\p(f)^*\p(f)\|<\e$ and 
$\|\p(ff^*)-\p(f)\p(f)^*\|<\e$ for $f\in F$. 
We say $\p$ is $(F,\e)$-suitable if it is 
$(F,\e)$-multiplicative and there are a faithful representation 
$C\subset\B(\hh)$ and a ucp map $\alpha\colon\M_p\to\B(\hh)$ 
such that $\|f-\alpha\p(f)\|<\e$ for $f\in F$. 
For a ucp map $\psi\colon C\to\M_q$, 
we denote by $\p\oplus\psi$ the ucp map 
$C\ni f\mapsto\diag(\p(f),\psi(f))\in\M_{p+q}$ 
and denote the $n$-fold sum $\p\oplus\cdots\oplus\p$ by $n\p$. 
For a ucp map $\psi\colon C\to\M_p$, we write 
$\p\simm_{F,\e}\psi$ if there is a unitary $u\in\M_p$ such that 
$\|\p(f)-u\psi(f)u^*\|<\e$ for $f\in F$. 
\end{defn}

Although, we are not going to use this fact, 
we remark that any ucp map $\p$ satisfies 
$$\|\p(ab)-\p(a)\p(b)\|\le\|\p(aa^*)-\p(a)\p(a)^*\|^{1/2}
\|\p(b^*b)-\p(b)^*\p(b)\|^{1/2}$$
for any $a$ and $b$ (cf.\ \cite{choi}). 
This justifies the term `$(F,\e)$-multiplicative'. 
The definition of $(F,\e)$-suitability is essentially same as 
the $(F,\e)$-admissibility of Dadarlat \cite{dsbq} and 
it was proved by Dadarlat \cite{djfa} that 
for any exact quasidiagonal $C^*$-algebra $C$, 
any finite subset $F\subset C$ and any $\e>0$, 
there is an $(F,\e)$-suitable ucp map from $C$ into a matrix algebra. 

We note the following easy facts as a lemma. 
\begin{lem}\label{note} 
If $\p\simm_{F,\e}\psi$ and $\psi\simm_{F,\delta}\rho$, then 
$\p\simm_{F,\e+\delta}\rho$. 
If $\p_j$ and $\psi_j$ are ucp maps such that 
$\p_j\simm_{F,\e}\psi_j$ for all $j=1,\ldots,k$, then 
$\bigoplus_{j=1}^k\p_j\simm_{F,\e}\bigoplus_{j=1}^k\psi_j$.
If $\p$ is $(F,\e)$-suitable and $\psi$ is $(F,\e)$-multiplicative, 
then $\p\oplus\psi$ is $(F,\e)$-suitable. 
\end{lem}

The following is our main tool. 
This is essentially due to Dadarlat \cite{dsbq}, 
but we include the proof for completeness. 
\begin{lem}\label{lem}
We put $h(\e)=5\sqrt{\e}$ and 
let $F$ be a finite set of unitaries in a unital $C^*$-algebra $C$ 
and $\e>0$. 
If $\p\colon C\to\M_p$ is an $(F,\e)$-suitable ucp map 
and $\psi\colon C\to\M_q$ is an $(F,\e)$-multiplicative ucp map, 
then there are $n\in\N$ and an $(F,h(\e))$-multiplicative 
ucp map $\rho\colon C\to\M_{np-q}$ such that 
$n\p\simm_{F,h(\e)}\psi\oplus\rho$. 
\end{lem}
\begin{proof} 
Let $C\subset\B(\hh)$ be a faithful representation 
and let $\alpha\colon\M_p\to\B(\hh)$ be a ucp map 
with $\| f-\alpha\p(f)\|<\e$ for $f\in F$. 
Let $\bar{\psi}\colon\B(\hh)\to\M_q$ be a ucp extension of $\psi$. 
Applying the Stinespring theorem to the ucp map 
$\bar{\psi}\alpha\colon\M_p\to\M_q$, we find $n\in\N$, 
an isometry $v\colon\ell_2^q\to\ell_2^n\otimes\ell_2^p$ 
such that 
$\bar{\psi}\alpha(x)=v^*(1\otimes x)v$ 
for all $x\in\M_p$. 
Thus denoting $\Phi=n\p$, we have that 
$\|v^*\Phi(f)v-\psi(f)\|<\e$ 
for $f\in F$. 
For the range projection $e=vv^*$ of $v$, it is routine to verify 
$\|(1-e)\Phi(f)e\|<2\sqrt{\e}$ and $\|e\Phi(f)(1-e)\|<2\sqrt{\e}$ 
for $f\in F$ (we recall that $f\in F$ is a unitary). 
Thus, defining a ucp map $\rho$ from $C$ into 
$\M_{np-q}=(1-e)(\M_n\otimes\M_p)(1-e)$ 
by $\rho(f)=(1-e)\Phi(f)(1-e)$, 
we are done. 
\end{proof} 

Given a unital separable exact $C^*$-algebra $A$, 
we define, for each $t\in[0,1]$, a unital exact contractible 
$C^*$-algebra $C(t)$ by $$C(t)=\{ f\in C([t,1],A) : f(1)\in\C1_A\}.$$
We note that $C:=C(0)$ is the unitization of the cone $CA$ over $A$ and 
$C(1)=\C$. 
We will prove that $C$ is AF-embeddable. 
It follows from Voiculescu's theorem \cite{voiculescu} that 
$C(t)$ is quasidiagonal for every $t\in[0,1]$. 
We denote by $\pi(t)\colon C\to C(t)$ the natural projection 
from $C$ onto $C(t)$ and denote by $\tilde{\sigma}(t)\colon C(t)\to C$ 
the $*$-homomorphism given by 
$$(\tilde{\sigma}(t)(f))(s)=\begin{cases}f(t) &\mbox{if }s<t\\ 
f(s) &\mbox{if }s\geq t\end{cases}$$
for $f\in C$. 
We note that $\sigma(t):=\tilde{\sigma}(t)\pi(t)$
is a continuous path of $*$-endomorphisms on $C$ 
with $\sigma(0)=\id_C$ and $\sigma(1)(C)=\C1_C$. 

\begin{lem}\label{key}
Let $C$ be as above. 
We give ourselves the following: 
\begin{enumerate} 
\item a finite set $F_0$ of unitaries in $C$, $\e_0>0$, 
\item $N_0\in\N$ so that 
$\|f(s)-f(t)\|_A<\e_0$ for any $f\in F_0$ and 
any $s,t\in[0,1]$ with $|s-t|\le\frac{1}{N_0}$, 
\item $\delta_0>0$ so that $h^{2N_0}(\delta_0)<\e_0$, 
(we recall that $h(\delta)=5\sqrt{\delta}$)
\item for each 
$t\in\{0,\frac{1}{N_0},\frac{2}{N_0},\ldots,\frac{N_0-1}{N_0}\}$, 
a $(\pi(t)(G_0),\delta_0)$-suitable ucp map $\tilde{\theta}_0(t)$ 
from $C(t)$ into a full matrix algebra $D_0(t)$, where 
$$G_0=\{\sigma(t)(f) : 
f\in F_0,\ t=0,\frac{1}{N_0},\frac{2}{N_0},\ldots,1\}\subset C$$ 
is the finite set of unitaries and $\tilde{\theta}_0(1):=\id_{C(1)}$. 
\end{enumerate}
Then, for any finite set $G_1$ of unitaries in $C$, 
$N_1\in\N$ which is a multiple of $N_0$ and $\delta_1>0$, 
there are, for each 
$t\in\{0,\frac{1}{N_1},\frac{2}{N_1},\ldots,\frac{N_1-1}{N_1}\}$,
a $(\pi(t)(G_1),\delta_1)$-suitable ucp map $\tilde{\theta}_1(t)$ 
from $C(t)$ into a full matrix algebra $D_1(t)$ 
and positive integers $n_t(s)\in\N$ 
for $s\in\{0,\frac{1}{N_0},\frac{2}{N_0},\ldots,1\}$ with $s\geq t$ 
such that 
$$\theta_1(t)\simm_{F_0,3\e_0}
\bigoplus_{s\geq t}n_t(s)\theta_0(s),$$
where $\theta_1(t)=\tilde{\theta}_1(t)\pi(t)$ and 
$\theta_0(s)=\tilde{\theta}_0(s)\pi(s)$ for every $s$.
\end{lem}
\begin{proof} 
We may assume that $G_1\supset G_0$ and $\delta_1<\delta_0$. 
For each 
$t\in\{0,\frac{1}{N_1},\frac{2}{N_1},\ldots,\frac{N_1-1}{N_1}\}$, 
we take a $(\pi(t)(G_1),\delta_1)$-suitable ucp map $\tilde{\nu}(t)$ 
from $C(t)$ into a full matrix algebra 
(which exists as it was shown by Dadarlat \cite{djfa}) 
and let $\nu(t)=\tilde{\nu}(t)\pi(t)$. 
We now fix 
$t\in\{0,\frac{1}{N_1},\frac{2}{N_1},\ldots,\frac{N_1-1}{N_1}\}$ 
and will construct $\theta_1(t)$. 

Let $k\in\N\cup\{0\}$ be such that $\frac{k-1}{N_0}<t\le\frac{k}{N_0}$. 
We note the crucial fact that $\sigma(\frac{j}{N_0})(G_0)\subset G_0$ 
for every $j=0,1,\ldots,N_0$. 
Applying Lemma \ref{lem} to $\tilde{\theta}_0(\frac{k}{N_0})$ and 
$\nu(t)\tilde{\sigma}(\frac{k}{N_0})$, 
we find $n_t(\frac{k}{N_0})$ and 
a $(\pi(\frac{k}{N_0})(G_0),h(\delta_0))$-multiplicative ucp map 
$\tilde{\rho}(\frac{k}{N_0})$ from $C(\frac{k}{N_0})$ 
into a full matrix algebra such that
$\rho(\frac{k}{N_0}):=\tilde{\rho}(\frac{k}{N_0})\pi(\frac{k}{N_0})$ 
satisfies 
$$\nu(t)\sigma(\frac{k}{N_0})\oplus\rho(\frac{k}{N_0})
\simm_{G_0,h(\delta_0)}
n_t(\frac{k}{N_0})\theta_0(\frac{k}{N_0}).$$
Applying Lemma \ref{lem} to $\tilde{\nu}(\frac{k}{N_0})$ and 
$\tilde{\rho}(\frac{k}{N_0})$, 
we find $m(\frac{k}{N_0})$ and 
a $(\pi(\frac{k}{N_0})(G_0),h^2(\delta_0))$-multiplicative ucp map 
$\tilde{\mu}(\frac{k}{N_0})$ from $C(\frac{k}{N_0})$ 
into a full matrix algebra such that 
$\mu(\frac{k}{N_0}):=\tilde{\mu}(\frac{k}{N_0})\pi(\frac{k}{N_0})$ 
satisfies
$$\rho(\frac{k}{N_0})\oplus\mu(\frac{k}{N_0})
\simm_{G_0,h^2(\delta_0)}
m(\frac{k}{N_0})\nu(\frac{k}{N_0}).$$
Since $\mu(\frac{k}{N_0})\tilde{\sigma}(\frac{k+1}{N_0})$ 
is $(\pi(\frac{k+1}{N_0})(G_0),h^2(\delta_0))$-multiplicative, 
we can further apply Lemma \ref{lem} and find 
$n_t(\frac{k+1}{N_0})$ and 
a $(\pi(\frac{k+1}{N_0})(G_0),h^3(\delta_0))$-multiplicative ucp map 
$\tilde{\rho}(\frac{k+1}{N_0})$ from $C(\frac{k+1}{N_0})$ 
into a full matrix algebra such that 
$\rho(\frac{k+1}{N_0})
:=\tilde{\rho}(\frac{k+1}{N_0})\pi(\frac{k+1}{N_0})$ 
satisfies 
$$\mu(\frac{k}{N_0})\sigma(\frac{k+1}{N_0})
\oplus\rho(\frac{k+1}{N_0})
\simm_{G_0,h^3(\delta_0)}
n_t(\frac{k+1}{N_0})\theta_0(\frac{k+1}{N_0}).$$

Continuing in this way, for each $j$ with $k<j<N_0$, 
we find $n_t(\frac{j}{N_0})$ and 
a $(G_0,h^{2(j-k)+1}(\delta_0))$-multiplicative 
ucp map $\tilde{\rho}(\frac{j}{N_0})$ from $C(\frac{j}{N_0})$ 
into a full matrix algebra such that 
$\rho(\frac{j}{N_0}):=\tilde{\rho}(\frac{j}{N_0})\pi(\frac{j}{N_0})$ 
satisfies
$$\mu(\frac{j-1}{N_0})\sigma(\frac{j}{N_0})
\oplus\rho(\frac{j}{N_0})
\simm_{G_0,h^{2(j-k)+1}(\delta_0)}
n_t(\frac{j}{N_0})\theta_0(\frac{j}{N_0})$$
and find $m(\frac{j}{N_0})$ and 
a $(G_0,h^{2(j-k+1)}(\delta_0))$-multiplicative 
ucp map $\tilde{\mu}(\frac{j}{N_0})$ from $C(\frac{j}{N_0})$ 
into a full matrix algebra such that 
$\mu(\frac{j}{N_0}):=\tilde{\mu}(\frac{j}{N_0})\pi(\frac{j}{N_0})$ 
satisfies 
$$\rho(\frac{j}{N_0})\oplus\mu(\frac{j}{N_0})
\simm_{G_0,h^{2(j-k+1)}(\delta_0)}
m(\frac{j}{N_0})\nu(\frac{j}{N_0}).$$
Finally, let $n_t(1)$ be so that 
$$\mu(\frac{N_0-1}{N_0})\sigma(1)=n_t(1)\theta_0(1)$$
(recall that $C(1)=\C$). 
Since 
$\|\sigma(\frac{j}{N_0})(f)-\sigma(\frac{j+1}{N_0})(f)\|_C<\e_0$ 
for every $f\in F_0$ and every $j=0,1,\ldots,N_0$, 
it follows from Lemma \ref{note} that 
\begin{align*}
\theta_1(t):=& 
\nu(t)\oplus m(\frac{k}{N_0})\nu(\frac{k}{N_0})
\oplus m(\frac{k+1}{N_0})\nu(\frac{k+1}{N_0})
\oplus\cdots\oplus
m(\frac{N_0-1}{N_0})\nu(\frac{N_0-1}{N_0})\\
\simm_{G_0,h^{2N_0}(\delta_0)}&
\nu(t)
\oplus\left(\rho(\frac{k}{N_0})\oplus\mu(\frac{k}{N_0})\right)
\oplus\left(\rho(\frac{k+1}{N_0})\oplus\mu(\frac{k+1}{N_0})\right)
\oplus\cdots\\ 
& \cdots\oplus
\left(\rho(\frac{N_0-1}{N_0})\oplus\mu(\frac{N_0-1}{N_0})\right)\\
\simm_{F_0,\e_0}&
\left(\nu(t)\sigma(\frac{k}{N_0})\oplus\rho(\frac{k}{N_0})\right)
\oplus
\left(\mu(\frac{k}{N_0})\sigma(\frac{k+1}{N_0})
\oplus\rho(\frac{k+1}{N_0})\right)
\oplus\cdots\\
& \cdots\oplus\left(\mu(\frac{N_0-2}{N_0})\sigma(\frac{N_0-1}{N_0})
\oplus\rho(\frac{N_0-1}{N_0})\right)
\oplus\mu(\frac{N_0-1}{N_0})\sigma(1)\\
\simm_{G_0,h^{2N_0}(\delta_0)}& 
n_t(\frac{k}{N_0})\theta_0(\frac{k}{N_0})
\oplus n_t(\frac{k+1}{N_0})\theta_0(\frac{k+1}{N_0})
\oplus\cdots\oplus n_t(\frac{N_0-1}{N_0})\theta_0(\frac{N_0-1}{N_0})
\oplus n_t(1)\theta_0(1).
\end{align*}
Since $\theta_1(t)$ is of the form 
$\theta_1(t)=\tilde{\theta}_1(t)\pi(t)$, 
we complete the proof by Lemma \ref{note}. 
\end{proof}

Now, it is not hard to show the proposition. 
\begin{proof}[Proof of Proposition \ref{prop1}] 
Fix a dense sequence $(f_i)_{i=1}^\infty$ 
in the set of unitaries of $C$ and let $F_k=\{ f_i : 1\le i\le k\}$.  
For each $k$, we take $N_k\in\N$ which is a multiple of $N_{k-1}$ 
such that 
$\|f(s)-f(t)\|_A<\frac{1}{2^k}$ for any $f\in F_k$ and 
any $s,t\in[0,1]$ with $|s-t|\le\frac{1}{N_k}$ 
and take $\delta_k>0$ so that $h^{2N_k}(\delta_k)<\frac{1}{2^k}$. 
Let 
$G_k=\{\sigma(t)(f) : 
f\in F_k,\ t=0,\frac{1}{N_k},\frac{2}{N_k},\ldots,1\}\subset C$ 
be a finite set of unitaries. 
By recursive use of Lemma \ref{key}, we find 
a $(\pi(t)(G_k),\delta_k)$-suitable ucp map 
$\tilde{\theta}_k(t)$ from $C(t)$ into a full matrix algebra $D_k(t)$ 
for each $k\in\N$ and 
$t\in\{0,\frac{1}{N_k},\frac{2}{N_k},\ldots,1\}$, 
with $\tilde{\theta}_k(1)=\id_{C(1)}$, 
and positive integers $n^{(k)}_t(s)$ for 
$t\in\{0,\frac{1}{N_{k+1}},\frac{2}{N_{k+1}},\ldots,1\}$ and 
$s\in\{0,\frac{1}{N_k},\frac{2}{N_k},\ldots,1\}$ with 
$s\geq t$ such that 
$$\theta_{k+1}(t)\simm_{F_k,\frac{3}{2^k}}
\bigoplus_{s\geq t}n^{(k)}_t(s)\theta_k(s),$$
where $\theta_k(t)=\tilde{\theta}_k(t)\pi_k(t)$ for every $k$ and $t$ 
(in particular, $n^{(k)}_1(1)=1$ for all $k$). 

We let $\theta_k=\bigoplus_t\theta_k(t)$ 
be the $(F_k,\frac{1}{2^k})$-suitable ucp map from $C$ into 
the finite dimensional $C^*$-algebra $D_k=\bigoplus_tD_k(t)$ 
(we apologize an abuse of the notation; 
$D_k$ is not a full matrix algebra). 
Then, the embedding $\iota_k$ of $D_k$ into $D_{k+1}$ 
is given along $(n^{(k)}_t(s))_{s,t}$ so that it satisfies 
$$\|\theta_{k+1}(f)-\iota_k\theta_k(f)\|<\frac{3}{2^k}$$
for $f\in F_k$. 
Now, one can embed $C$ into the AF-algebra 
$D=\overline{\bigcup_{k\in\N}D_k}$ through $(\theta_k)_{k\in\N}$. 
(This embedding is indeed multiplicative as all unitaries in $C$ are 
contained in the multiplicative domain \cite{choi}.) 
This completes the proof. 
\end{proof} 

\section{Proof of Proposition \ref{prop2}}
The proof goes as that of Theorem 3.9 
in Spielberg \cite{spielberg}, but we avoid the UCT. 
\begin{proof}[Proof of Proposition \ref{prop2}]
Let $\p\colon A\to B$ and $\psi\colon B\to A$ be $*$-homomorphisms 
such that there is a continuous path $(\gamma(t))_{t\in[0,1]}$ 
of $*$-endomorphisms on $A$ 
with $\gamma(0)=\id_A$ and $\gamma(1)=\psi\p$. 
Let $CA=C_0([0,1),A)$ be the cone over $A$ 
and let $Z_\psi=\{ (f,b)\in C([0,1],A)\oplus B : f(1)=\psi(b)\}$ 
be the mapping cylinder of $\psi$. 
Then, we have a short exact sequence 
$0\to CA\to Z_\psi\to B\to0$
which clearly splits. 
We note that $A$ is embedded into $Z_\psi$ through $\gamma\oplus\p$. 
Since $CA$ is AF-embeddable, we may find 
the following commutative diagram (cf.\ 1.11 in \cite{spielberg}) 
with exact rows: 
$$\begin{CD}
0 @>>> CA @>>> Z_\psi @>>> B @>>> 0\\
@. @VVV @VVV @| @.\\
0 @>>> J @>>> E @>>> B @>>> 0
\end{CD}$$
where the vertical $*$-homomorphisms are all injective 
and $J$ is an AF-algebra.  
Since the lower row again splits, we obtain the AF-embeddability 
of $E$ as in Lemma 1.13 in \cite{spielberg}. 
\end{proof} 
\providecommand{\bysame}{\leavevmode\hbox to3em{\hrulefill}\thinspace} 

\end{document}